\documentclass[12pt]{amsart}

\usepackage{amssymb}
\usepackage{amsfonts,cmtiup,comment}

\makeatletter
\def\blfootnote{\xdef\@thefnmark{}\@footnotetext}
\makeatother

\newtheorem{theorem}{Theorem}[section]
\newtheorem{lemma}[theorem]{Lemma}

\newtheorem{corollary}[theorem]{Corollary}

\theoremstyle{remark}
\newtheorem{remark}[theorem]{Remark}

\newcommand{\bt}{\begin{theorem}}
\newcommand{\et}{\end{theorem}}
\newcommand{\be}{\begin{equation}}
\newcommand{\ee}{\end{equation}}
\newcommand{\bc}{\begin{corollary}}
\newcommand{\ec}{\end{corollary}}
\newcommand{\bp}{\begin{proof}}
\newcommand{\ep}{\end{proof}}
\newcommand{\bl}{\begin{lemma}}
\newcommand{\el}{\end{lemma}}
\newcommand{\br}{\begin{remark}}
\newcommand{\er}{\end{remark}}

\let\leq=\leqslant
\let\geq=\geqslant
\newcommand{\ed}{\end{document}}

\begin{document}

\title{On the length of finite factorized groups}

\author{E. I. Khukhro and P. Shumyatsky}

\keywords{Finite groups, nonsoluble length, generalized Fitting height, factorized group}
\subjclass{20D40}

\begin{abstract} 
The nonsoluble length~$\lambda (G)$ of a finite group~$G$ is defined as the number of nonsoluble factors in a shortest normal series
each of whose factors either is soluble or is a direct product of nonabelian simple groups. The generalized Fitting height of a finite group~$G$ is
the least number $h=h^*(G)$ such that $F^*_h(G)=G$, where  $F^*_1(G)=F^*(G)$ is the generalized Fitting subgroup, and $F^*_{i+1}(G)$ is the inverse image of $F^*(G/F^*_{i}(G))$. It is proved that if a finite group $G=AB$ is factorized by two subgroups of coprime orders, then the nonsoluble length of~$G$ is bounded in terms of the generalized Fitting heights of~$A$ and~$B$. It is also proved that if, say, $B$ is soluble of derived length~$d$, then the generalized Fitting height of~$G$ is bounded in terms of~$d$ and the generalized Fitting height of~$A$.
\end{abstract}

\maketitle

\blfootnote{This work was supported by CNPq-Brazil. The first author thanks CNPq-Brazil and the University of Brasilia for support and hospitality that he enjoyed during his visits to Brasilia.}

\section{Introduction}

If a finite group $G=AB$ is a product of two subgroups $A,B$, then it is natural to expect restrictions on the structure of~$G$ in terms of~$A$ and~$B$. There are many papers dedicated to such restrictions, like criteria for being non-simple, or bounds for certain parameters when $G$ is soluble. We refer the reader to the books \cite{AFD} and \cite{BB} and extensive bibliography therein.

In
the present paper we consider finite factorizable groups $G=AB$ in relation to two important parameters that in a sense measure the complexity of nonsoluble groups. One of these parameters is based on the generalized Fitting subgroup~$F^*(G)$. Recall that $F^*(G)$ is the product of the Fitting subgroup~$F(G)$ and all subnormal quasisimple subgroups; here a group is quasisimple if it is perfect and its quotient by the centre is a non-abelian simple group. Then the \textit{generalized Fitting series} of~$G$ is defined starting from $F^*_1(G)=F^*(G)$, and then by induction, $F^*_{i+1}(G)$ being the inverse image of $F^*(G/F^*_{i}(G))$. The least number~$h$ such that $F^*_h(G)=G$ is naturally defined as the \textit{generalized Fitting height}~$h^*(G) $ of~$G$. Clearly, if $G$ is soluble, then $h^*(G)=h(G)$ is the ordinary Fitting height of~$G$. Bounding the generalized Fitting height of a finite group~$G$ greatly facilitates using the classification of finite simple groups (and is itself often obtained using the classification). Examples include reduction theorems for problems of Restricted Burnside type  \cite{ha-hi, wil83}, which included in effect bounding the generalized Fitting height.

Another length parameter is
the \textit{nonsoluble length}~$\lambda (G)$ of a finite group~$G$, which is defined as the number of nonsoluble factors in a shortest normal series
each of whose factors either is soluble or is a direct product of nonabelian simple groups. Bounding the nonsoluble length implicitly appeared in \cite{ha-hi,wil83} and, more recently, in \cite{khu-shu131,khu-shu132,khu-shu133}.

Our main result is a bound for the nonsoluble length of a finite group factorized by two subgroups of coprime orders in terms of the generalized Fitting height of the factors.

\bt\label{t1}
Suppose that a finite group $G=AB$ admits a factorization by two subgroups $A$,~$B$ of coprime orders. Then the nonsoluble length~$\lambda (G)$ of~$G$ is bounded in terms of the generalized Fitting heights $h^*(A)$ and $h^*(B)$ of the factors. More precisely, $\lambda (G)\leq 2^{h^*(A)+h^*(B)}-1$.
\et

The proof of Theorem~\ref{t1} uses Schreier's conjecture on solubility of outer automorphism groups of non-abelian finite simple groups, which was confirmed by the classification. Examples show that in Theorem~\ref{t1} the generalized Fitting height of~$G$ cannot be bounded in terms of the generalized Fitting heights of the factors. There are also examples showing that the nonsoluble length cannot be bounded in terms of the nonsoluble length of the factors.

Since one of the subgroups of coprime orders is soluble by the Feit--Thompson theorem, it makes sense to consider the derived length of a soluble factor. Then an estimate for the generalized Fitting height of the group $G=AB$ can be obtained.

\bc\label{t2}
Suppose that a finite group $G=AB$ admits a factorization by two subgroups $A$,~$B$ of coprime orders, of which $B$ is soluble of derived length~$d$. Then the generalized Fitting height $h^*(G)$ of~$G$ is bounded in terms of~$d$ and  $h^*(A)$. \ec

Corollary~\ref{t2} follows from Theorem~\ref{t1} and the result of
Casolo, Jabara, and Spiga \cite{cajasp}, which gives a bound for the Fitting height when $G$ is soluble.

It is not clear at the moment how essential the coprimeness conditions are in Theorem~\ref{t1} and Corollary~\ref{t2}. Note that by Kegel's theorem  \cite{keg61} a finite group factorized by two nilpotent subgroups of not necessarily coprime orders is soluble.

We conclude the Introduction by examples mentioned above. It is clear that the (generalized) Fitting height of $G=AB$ cannot be bounded in terms of the (generalized) Fitting heights  of the factors, since there are groups of order $p^aq^b$ of arbitrarily large Fitting height. In order to show that the nonsoluble length of  $G=AB$ cannot be bounded in terms of the nonsoluble lengths of the factors, consider the alternating group $A_5$, which has Hall $\{2,3\}$- and $5$-subgroups. Then a repeated wreath product $G=A_5\wr (A_5\wr \dots \wr(A_5\wr A_5)\dots ))$ also has Hall $\{2,3\}$- and $5$-subgroups, say, $A$ and $B$, which are soluble. Then $G=AB$. But the nonsoluble length of $G$ is unbounded.

\section{Preliminaries}
\label{s-prel}

Recall that the generalized Fitting subgroup $F^*(G)$ of a finite group~$G$ is the product of the Fitting subgroup $F(G)$ and the characteristic subgroup $E(G)$, which is a central product of all subnormal quasisimple subgroups of~$G$ (components of~$F^*(G)$), that is, $E(G)=\prod Q_i$ over all $Q_i$ such that $Q_i$ is subnormal in~$G$, $Z(Q_i)\leq [Q_i,Q_i]$, and $S_i=Q_i/Z(Q_i)$ is a non-abelian simple group. Then $[F(G),E]=1$ and $E(G)/Z(E(G))\cong F^*(G)/F(G)$ is the direct product of the~$S_i$. Acting by conjugation, the group~$G$ permutes the factors~$Q_i$ and $C_G(F^*(G))\leq F(G)$. The following fact (see, for example, \cite[Lemma~2.1]{khu-shu131}) is a well-known consequence of Schreier's conjecture on solubility of outer automorphism groups of non-abelian finite simple groups confirmed by the classification.

\bl\label{l21}
Let $L/S(G)=F^*(G/S(G))$ be the generalized Fitting subgroup of the quotient by the soluble radical $S(G)$ of a finite group~$G$, and let $K$ be the kernel of the permutational action of~$G$ on the set of subnormal simple factors of~$L/S(G)$. Then $K/L$ is soluble. \qed
\el

We shall use without special references the following well-known properties of the generalized Fitting subgroups relative to normal subgroups: if $N\unlhd {}G$, then $F^*(N)=F^*(G)\cap N$ and $F^*(G)N/N \leq F^*(G/N)$.
These and other properties follow, for example, from the fact that $F^*(G)$ is the set of all elements of~$G$ that induce inner automorphisms on every chief factor of~$G$; see, for example, \cite[Ch.~X, \S\,13]{hup3}. It is easy to see that similar properties hold for the higher terms of the generalized Fitting series: if $N\unlhd G$, then $F^*_i(N)=F^*_i(G)\cap N$ and $F^*_i(G)N/N \leq F^*_i(G/N)$. Hence, $h^* (N)\leq h^* (G)$, $h^* (G/N) \leq h^* (G)$, and $h^* (G)\leq h^* (N)+h^* (G/N)$.

We shall need the following technical lemma.

\begin{lemma}\label{lf*}
Suppose that $T$ is a subgroup of a finite group~$H$ such that its normal closure $\langle T^H\rangle$ is equal to the direct product of its conjugates $\langle T^H\rangle=T^{x_1}\times \dots \times T^{x_n}$, all its conjugates appear as these factors, that is, $T^y\in\{T^{x_1}, \dots , T^{x_n}\}$ for any $y\in H$, and $H=F^*(H)\langle T^H\rangle$. Then either $T$ is normal in~$H$, or $T$ is a $p$-group for some prime~$p$.
\end{lemma}

\bp
Suppose that $T$ is nonsoluble. Consider the quotient $\bar H= H/S(\langle T^H\rangle)$ by the soluble radical of $\langle T^H\rangle$, which is of course equal to the product $S(T)^{x_1}\times \dots \times S(T)^{x_n}$. Since $\overline{F^*(H)}\leq F^*(\bar H)$, clearly, $\bar H=F^*(\bar H)\langle \bar T^{\bar H}\rangle$, where $1\ne \bar T=T/S(T)$ and
$F(\bar T)=1$. Then the quasisimple components of~$F^*(\bar T)$ are also components of~$F^*(\bar H)$ because $\bar T$ is subnormal in~$\bar H$. These components are normal in~$F^*(\bar H)$, and therefore $\bar T$ is normalized by~$F^*(\bar H)$, so that $TS(\langle T^H\rangle )$ is normalized by~$F^*( H)$. Since  $T^y\in\{T^{x_1}, \dots , T^{x_n}\}$ for any $y\in H$ by hypothesis and $T$ is nonsoluble by assumption, it follows that  $T$ itself is  normalized by~$F^*( H)$ and therefore is normal in~$H$.

Thus we can assume that $T$ is soluble and $O_p(T)\ne 1$ for some prime~$p$. Then $O_p(T)^{x_1}\times \dots \times O_p(T)^{x_n}= O_p(\langle T^H\rangle)\leq O_p(H)$. Since all $p'$-elements of~$F^*(H)$ centralize~$O_p(H)$, they all normalize each factor~$T^{x_i}$. Suppose that $T\ne O_p(T)$. Then in the quotient $\tilde H= H/O_p(\langle T^H\rangle)$ we shall have $O_q(\tilde T)\ne 1$ for a different prime~$q$. For the same reasons, $p$-elements of the image of~$F^*(H)$, which is contained in $F^*(\tilde H)$, will be $q'$-elements centralizing $O_q(\tilde T)$ and therefore normalizing~$\tilde T$. As a result, the whole $F^*(H)$  normalizes~$TO_p(\langle T^H\rangle)$.  Since  $T^y\in\{T^{x_1}, \dots , T^{x_n}\}$ for any $y\in H$ by hypothesis and $T$ is not a $p$-group by assumption, then $F^*(H)$ normalizes~$T$ and therefore $T$ is normal in~$H$.
\ep

Recall that the nonsoluble length~$\lambda (G)$ of a finite group~$G$ is the number of nonsoluble factors in a shortest normal series
each of whose factors either is soluble or is a direct product of nonabelian simple groups. For example, $\lambda (K)\leq 1$ for the subgroup~$K$ in Lemma~\ref{l21}. This parameter also behaves well in relation to normal subgroups: if $N\unlhd {}G$, then $\lambda (N)\leq \lambda (G)$, $\lambda (G/N) \leq \lambda G$, and $\lambda (G)\leq \lambda (N)+\lambda (G/N)$.

The following two lemmas on factorizable groups must be well-known.

\bl\label{l1}
 Let $a$, $b$ be elements of coprime orders of a finite group~$G$ such that their
 product $ab$ belongs to a normal subgroup~$N$. Then both $a$ and~$b$ belong to~$N$.
\el

\bp
Indeed, the images of $a,b^{-1}$ in~$G/N$ are equal and at the same time have coprime orders. Hence both are trivial.
\ep

\bl\label{l11}
Let $G=AB$ be a finite group factorized by subgroups $A,B$ of coprime order, and let $N$ be a normal subgroup of~$G$. Then $N=(N\cap A)(N\cap B)$.
\el

\bp
Every element $n\in N$ is a product $n=ab$ for some $a\in A$, $b\in B$. By Lemma~\ref{l1}, both $a\in N$ and $b\in N$.
\ep

We shall use the following theorem of Casolo, Jabara, and Spiga.

\begin{theorem}[{{\rm \cite[Theorem~1.1]{cajasp}}}]\label{t3}
Let $G=AB$ be a finite soluble group factorised by its proper subgroups~$A$ and~$B$ with $\gcd(|A|,|B|)=1$. Then $h(G)\leq h(A)+h(B)+4d(B)-1$.
Moreover, if $|B|$ is odd, then $h(G)\leq h(A)+h(B)+2d(B)-1$,
and if $B$ is nilpotent, then $h(G)\leq h(A)+2d(B)$.
\end{theorem}

\section{Proofs of the main results}
\label{s-main}

\bp[Proof of Theorem~\ref{t1}] We have a finite group $G=AB$ factorized by subgroups $A,B$ of coprime orders; we wish to obtain a bound for the nonsoluble length $\lambda (G)$ of~$G$ in terms of the generalized Fitting heights~$h^*(A)$ and~$h^*(B)$ of~$A$ and~$B$. We are actually going to prove that $\lambda (G)\leq 2^{h^*(A)+h^*(B)}-1$. We use induction on $h^*(A)+h^*(B) $. When $h^*(A)+h^*(B) =1$, one of the groups $A,B$ is trivial, so the result is obviously true. In the general case we can of course assume that $S(G)=1$. Then $F^*(G)=S_1\times \dots \times S_m$ for non-abelian simple groups~$S_i$, which are permuted by~$G$ acting by conjugation.

By Lemma~\ref{l11} we have $F^*(G)=(F^*(G)\cap A)(F^*(G)\cap B)$. By applying Lemma~\ref{l11} again we get $S_i=(S_i\cap A)(S_i\cap B)$ for every $i=1,\dots ,m$. We denote for short $S_{iA}=S_i\cap A$ and $S_{iB}=S_i\cap B$.

\bl\label{l-perm}
Acting by conjugation, the subgroup~$A$ permutes the~$S_{iA}$, and $B$ permutes the~$S_{iB}$.
\el

\bp For $a\in A$ we have $S_{i}^a=S_j$ for some~$j$ since $G$ permutes the~$S_i$. Then $S_{iA}^a=(S_i\cap A)^a=S_j\cap A=S_{jA}$. The argument for~$B$ and the~$S_{iB}$ is similar.
\ep

\bl\label{l-order}
If $S_i^g=S_j$ for some $g\in G$, then $|S_{iA}|=|S_{jA}|$ and $|S_{iB}|=|S_{jB}|$.
\el

\bp Let $g=ab$ for $a\in A$, $b\in B$. Let $S_i^a=S_k$. Then $S_{iA}^a=S_{kA}$ by Lemma~\ref{l-perm}. Hence, $|S_{iB}|=|S_i:S_{iA}|=|S_k:S_{kA}|=|S_{kB}|$. A similar calculation for $S_k^b=S_j$ shows that  $|S_{kB}|=|S_{jB}|$ and  $|S_{kA}|=|S_{jA}|$.
\ep

A key to the induction step in the proof of Theorem~\ref{t1} is given by the following lemma.

\bl\label{l2}
If $S_{iA}$ is not a $p$-group (for any prime~$p$), then $S_i$ is normalized by~$F^*(A)$, and if  $S_{iB}$ is not a $p$-group (for any prime~$p$), then $S_i$ is normalized by~$F^*(B)$.
\el

\bp
We prove the lemma for~$A$, since the proof is the same for~$B$. We argue by contradiction: suppose that $F^*(A)$ does not normalize~$S_{i}$. Then $F^*(A)$ also does not normalize~$S_{iA}$. Then $H=\langle F^*(A), S_{iA}\rangle =F^*(A)(S_{iA}^{x_1}\times\dots\times S_{iA}^{x_k})$ for $k>1$ and for some elements $x_s\in F^*(A)$. This subgroup~$H$ satisfies the hypotheses of Lemma~\ref{lf*}, by which $S_{iA}$ is $p$-group for some prime~$p$, a contradiction.
\ep

We now complete the proof of Theorem~\ref{t1}. Let
$$K_A=\bigcap \{N_G(S_i)\mid S_{iA}\text{ is not a $p$-group}\} $$
be the intersection of the normalizers of all the~$S_i$ such that $S_{iA}$ is not a $p$-group (for any~$p$).
Similarly, let
$$K_B=\bigcap \{N_G(S_i)\mid S_{jB}\text{ is not a $p$-group}\} $$
be the intersection of the normalizers of all the~$S_j$ such that $S_{jB}$ is not a $p$-group (for any~$p$).
By Lemma~\ref{l2}, we have $F^*(A)\leq K_A$ and $F^*(B)\leq K_B$.
It follows from Lemmas~\ref{l-perm} and~\ref{l-order} that both $K_A$ and~$K_B$ are normal subgroups of~$G$. Let
$$K=\bigcap _{i=1}^mN_G(S_i),$$
which is a normal subgroup of nonsoluble length at most 1 by Lemma~\ref{l21}.
Since every~$S_i$ is a non-abelian simple group, it cannot be a product of two Sylow subgroups by the Burnside $p^aq^b$ theorem. Therefore,  for each~$i$, at least one of~$S_{iA}$ or $S_{iB}$ is not a $p$-group. Hence, $K_A\cap K_B=K$.

The quotient $\bar G=G/K_A=\bar A\bar B$ satisfies the hypothesis of the theorem with $h^*(\bar A)+h^*(\bar B)\leq h^*(A)+h^*( B)-1$, since $F^*(A)\leq K_A$. By induction, $\lambda (\bar G)\leq 2^{h^*(A)+h^*( B)-1}-1$. Similarly, $\lambda (G/K_B)\leq 2^{h^*(A)+h^*( B)-1}-1$. Since the quotient $K_A/K=K_A/(K_A\cap K_B)$ is isomorphic to the normal subgroup $K_AK_B/K_B$ of~$G/K_B$, we also have $\lambda (K_A/K)\leq 2^{h^*(A)+h^*( B)-1}-1$. Taking the sum over the normal series $1\leq K\leq K_A\leq G$, we obtain the required bound $\lambda (G)\leq 2^{h^*(A)+h^*( B)-1}-1+ 2^{h^*(A)+h^*( B)-1}-1+1= 2^{h^*(A)+h^*( B)}-1$.
\ep

\bp[Proof of Corollary~\ref{t2}] We have a finite group $G=AB$ factorized by two subgroups $A$,~$B$ of coprime orders, of which $B$ is soluble of derived length~$d$. We wish to bound the generalized Fitting height~$h^*(G)$ of~$G$ in terms of~$d$ and~$h^*(A)$. By Theorem~\ref{t1} the group~$G$ has a normal series of length bounded in terms of~$h^*(A)$ and~$d$ each factor of which either is soluble or is a direct product of non-abelian simple groups.
Since the generalized Fitting height of a direct product of simple groups is 1, it remains to bound the Fitting height of the soluble normal sections. Let $M\geq N$ be normal subgroups of~$G$ with $M/N$ soluble. By Lemma~\ref{l11} we have $M=(M\cap A)(M\cap B)$, so that $M/N=((M\cap A)N/N)((M\cap B)N/N)$. The quotients $(M\cap A)N/N$ and $(M\cap B)N/N$ are soluble. Therefore the Fitting height of the first of them is equal to the generalized Fitting height, which is at most~$h^*(A)$. The derived length of the second is at most~$d$.
It remains to apply Theorem~\ref{t3} of Casolo, Jabara, and Spiga \cite{cajasp}.
\ep

\end{document}